\newif\ifreport
\newif\ifnotreport
\reporttrue 
\ifreport
	\notreportfalse
\else
	\notreporttrue
\fi

\ifreport
	\documentclass{article}
\else
	\documentclass{scspaperproc}
\fi 

\usepackage{latexsym}
\usepackage{mathptmx}

\usepackage{fixltx2e}
\usepackage[T1]{fontenc}
\usepackage{ae,aecompl}
\usepackage[utf8]{inputenc}
\usepackage[english]{babel}
\usepackage{verbatim}
\usepackage{graphicx}
\usepackage{amsfonts}
\usepackage{amsmath}
\usepackage{amssymb}
\usepackage{wrapfig}
\usepackage{stmaryrd}
\usepackage{algorithm}
\usepackage{algpseudocode}
\usepackage{todonotes}
\usepackage{xcolor}
\usepackage{paralist}
\usepackage[pdftex,colorlinks=true,urlcolor=blue,citecolor=black,anchorcolor=black,linkcolor=black]{hyperref}
\usepackage[nameinlink,capitalise]{cleveref}

\usepackage{hyphenat}
\graphicspath{{imgs/}}

\newcommand{\algoref}[1]{Algorithm~\ref{#1}}

\newcommand{\concept}[2]{%
  \phantomsection
  #1\def\@currentlabel{\unexpanded{#1}}\label{#2}%
}

\newcommand{\footurl}[1]{\footnote{\url{#1}}}

\newcommand{\brackets}[1]{\ensuremath{ \left[ #1 \right] }}
\newcommand{\tuple}[1]{\ensuremath{ \left\langle #1 \right\rangle }}
\newcommand{\set}[1]{\ensuremath{ \left\{ #1 \right\}}}

\newcommand{\pargroup}[1]{\ensuremath{ \left( #1 \right)}}
\newcommand{\dert}[1]{\ensuremath{#1'}}
\newcommand{\dertb}[1]{\ensuremath{\brackets{#1}'}}
\newcommand{\dernt}[2]{\ensuremath{#1^{(#2)}}}
\newcommand{\inftyint}[1]{\ensuremath{\int_{-\infty}^{\infty}#1}}

\newcommand{\setreal}[0]{\ensuremath{\mathbb{R}}}

\newcommand{\floor}[0]{\ensuremath{\mathit{floor}}}

\newcommand{\sumseq}[1]{\sum_{#1 \in \set{#1}}}
\newcommand{\conditional}[3]{%
\begin{cases}
	#1 \text { if } #3 \\
	#2 \text { otherwise }
\end{cases}}

\newenvironment{aligneq}%
{
\begin{equation}
\begin{aligned}
}{
\end{aligned}
\end{equation}
}

\begin{document}

\ifnotreport
	\setlength{\belowdisplayskip}{5pt}
	\setlength{\abovedisplayskip}{5pt}
\fi

\ifnotreport
	\SCSpagesetup{Gomes, Van Tendeloo, Denil, De Meulenaere, and Vangheluwe}

	\def\SCSconferenceacro{SpringSim}
	\def\SCSpublicationyear{2017}
	\def\SCSconferencedates{April 23-26}
	\def\SCSconferencevenue{Virginia Beach, VA, USA}
	\def\SCSsymposiumacro{TMS/DEVS}
\fi

\title{Hybrid System Modelling and Simulation with Dirac Deltas}

\author{Cláudio Gomes \\
University of Antwerp \\
Claudio.Gomes@uantwerp.be \\
\and Yentl Van Tendeloo \\
University of Antwerp \\
Yentl.VanTendeloo@uantwerp.be\\
\and Joachim Denil \\
University of Antwerp \\
Flanders Make \\
Joachim.Denil@uantwerp.be \\
\and Paul De Meulenaere \\
University of Antwerp \\
Flanders Make \\
Paul.DeMeulenaere@uantwerp.be \\
\and Hans Vangheluwe \\
University of Antwerp \\
McGill University \\
Flanders Make \\
Hans.Vangheluwe@uantwerp.be
}


\maketitle

\ifreport
	\tableofcontents
\fi

\ifnotreport
	\vspace{-30pt}
\fi

\ifreport
	\newcommand{\shortcite}[1]{\cite{#1}}
\fi

\section*{Abstract}

For a wide variety of problems, creating detailed continuous models of (continuous) physical systems is, at the very least, impractical.  
Hybrid models can abstract away short transient behaviour (thus introducing discontinuities) in order to simplify the study of such systems.
For example, when modelling a bouncing ball, the bounce can be abstracted as a discontinuous change of the velocity, instead of resorting to the physics of the ball (de-)compression to keep the velocity signal continuous.
Impulsive differential equations can be used to model and simulate hybrid systems such as the bouncing ball.
In this approach, the force acted on the ball by the floor is abstracted as an infinitely large function in an infinitely small interval of time, that is, an impulse.
Current simulators cannot handle such approximations well due to the limitations of machine precision.

In this paper, we explore the simulation of impulsive differential equations, where impulses are first class citizens.
We present two approaches for the simulation of impulses: symbolic and numerical.
Our contribution is a theoretically founded description of the implementation of both approaches in a Causal Block Diagram modelling and simulation tool. 
Furthermore, we investigate the conditions for which one approach is better than the other.

\textbf{Keywords:} Dirac delta, hybrid system, distribution theory, causal block diagrams.

\section{Introduction}

A model of a physical system\footnote{In this work, we consider physical systems to be continuous. This excludes Quantum systems.} has to pragmatically capture only the essential features for a task at hand \shortcite{Cellier1991}.
If the task is to analyze some behavior at a macro time scale, it is impractical -- or even currently impossible \shortcite{Stewart2000} -- to create high fidelity continuous models that also explain the behavior at micro time scales.
Abstraction is the only way to deal with complexity and create useful models \shortcite{Mosterman1998a}. 

As a simple example, consider a ball with unit mass that bounces once on the floor. 
The macro time scale behavior is the overall trajectory of the ball, whereas the micro time scale behavior captures the compression of the ball in the floor.
It can be modelled by means of first order Ordinary Differential Equations (ODEs):
\begin{aligneq}\label{eq:bouncing_ball_unabstracted}
&\dernt{y}{2} = -g + F_c(t)
 \ \ \text{with} \ \  
y(0) = y_0
 \ \ \text{and} \ \  
\dert{y}(0) = v_0
\end{aligneq}
where $\dernt{y}{2}$ denotes the acceleration of the ball (second derivative of the position), and $y_0$ and $v_0$ are the initial conditions of the dynamical model (constants). $F_c(t)$ denotes the force impinged by the floor as a result of the collision, which we want to abstract.
There are two possible approaches to find the solution $y(t)$ that satisfies \cref{eq:bouncing_ball_unabstracted}, both based on using conservation of momentum to ensure the abstraction of $F_c(t)$ is correct.

\textbf{Separation of dynamics} separates the time continuum into two groups: before and after the contact. Momentum conservation gives the velocity of the ball immediately after the collision. Let $t_c$ denote that time at which the contact occurs. \cref{eq:bouncing_ball_unabstracted} is then split in two:
\begin{aligneq}\label{eq:bouncing_ball_model_separated}
	\begin{pmatrix}
	\dernt{y}{2}	&= -g \\
	y(0)			&= y_0 \\ 
	\dert{y}(0)		&= v_0 
	\end{pmatrix}
	\ \ \ \text{for $0 \leq t < t_c$, and} \ \ 
	\begin{pmatrix}
	\dernt{y}{2}	&= -g \\
	y(t_c)			&= y(t_c^-) \\
	\dert{y}(t_c)	&= -\dert{y}(t_c^-)
	\end{pmatrix}
	\ \ \ \text{  for $t \geq t_c$}
\end{aligneq}
These are solved separately. $t_c$ can be found from the solution to the leftmost ODE.
A simulator can be used to compute the solution in the interval $t \in [0, t_c)$. 
Afterward, the same simulator can compute the solution in the interval $0 \leq t < t_c$, with the new initial conditions.

The activity of re-setting the same simulator with new initial conditions is called re-initialization and it is a common approach to the simulation of systems with discontinuities \shortcite{Cellier2006}.

The \textbf{Impulsive forces} approach acknowledges that, since the velocity is the integral of the acceleration, during the collision, $F_c(t)$ must be large enough to cause complete ``inversion'' of the velocity, even if its shape is unknown.
To be concrete, let $t_c^-$ denote the time at which the collision starts and $t_c^+$ the time at which the collision ends. 
Then, after the collision, the velocity of the ball is given by
\begin{aligneq}\label{eq:colision_before_after}
\dert{y}(t_c^+) = \dert{y}(t_c^-) + \int_{t_c^-}^{t_c^+} -g + F_c(\tau) d\tau
\end{aligneq}
Conservation of momentum implies that $\dert{y}(t_c^+) = - \dert{y}(t_c^-)$, and the effect that $-g$ has on the integral can be neglected. Thus:
\begin{aligneq}\label{eq:integral_requirements}
\int_{t_c^-}^{t_c^+} F_c(\tau) d\tau = - 2 \dert{y}(t_c^-)
\end{aligneq}
\emph{independently of the shape} that $F_c(\tau)$ assumes.
The Dirac impulse \shortcite{Dirac1981} formalizes the abstraction of a large continuous function, whose integral is a known value over a small interval of time.
More details will be given in \cref{sec:background} but it suffices to postulate that $\int_{0^-}^{0^+} \delta(\tau) d\tau = 1$.
Hence, setting $F_c(\tau) = - 2 \dert{y}(t_c^-) \delta(t - t_c)$ satisfies \cref{eq:integral_requirements} and  \cref{eq:colision_before_after}. 
With this approach, \cref{eq:bouncing_ball_unabstracted} is rewritten to
\begin{aligneq}\label{eq:bouncing_ball_dirac}
&\dernt{y}{2} = -g - 2 \dert{y}(t_c^-) \delta(t - t_c) 
 \ \ \text{with} \ \  
y(0) = y_0
 \ \ \text{and} \ \  
\dert{y}(0) = v_0
\end{aligneq}
which is an impulsive differential equation \shortcite{Lakshmikantham1989}. 

In the bouncing ball example, both the separation of dynamics and impulse based approaches assume that the duration of the contact is infinitely small and are able to reach the same solution.
The separation of dynamics allows the modeller more freedom in setting the initial states after an impulse occurs (e.g., the initial condition in \cref{eq:bouncing_ball_model_separated} could have been $y(t_c) = 0$ instead of $y(t_c) = y(t_c^-)$).
This means more mistakes can be made (e.g., violating physical laws).

The re-initialization technique can be used for the simulation of impulsive differential equations.
However, for complex models where impulses are regularly exchanged, a direct manipulation -- symbolic or numerical -- of impulses avoids the need to reset the whole state of the simulator.
A symbolic manipulation means that impulses are encoded as first class citizens in the operations of the simulator.
And by numerical we mean that impulses are approximated as very large values.

Simulators with support for Dirac impulses have been proposed in \shortcite{Nilsson2003} and more recently in \shortcite{Lee2015}, but both approaches deal with the symbolic computation of impulses, and not their numerical approximation, nor evaluation or comparison is advanced for the approach.

Our \emph{research hypothesis} is: (RH) a simulator which manipulates Dirac deltas symbolically (such as the one proposed in \shortcite{Nilsson2003}) is more accurate than one that just operates with approximated impulses.

Our \emph{contribution} can be divided into: 
(C1) a derivation of the decision and inversion symbolic operations performed on Dirac impulses;
(C2) a numerical approximation of impulses;
(C3) the conditions for which manipulating impulses symbolically yields more accurate solutions than the numerical approximations;

C1 builds on the work of \shortcite{Nilsson2003}, where many of the symbolic operations on Dirac deltas are derived.
Here, we derive the remaining operations -- decision and inversion -- and provide justification to the use of distributions as a formalization of a correct impulse symbolic simulator.
C2 and C3 answer RH.
An implementation of a simulator capable of simulating Causal Block Diagrams with impulses, symbolically and numerically, as well as the experiments used in this paper, are available at 
\url{http://msdl.cs.mcgill.ca/people/claudio/diraccbds}.

The next section provides some background and \cref{sec:hybridcbds:semantics} derives the main symbolic operations on impulses.
\cref{sec:general_cbds_discussion} describes how impulses are approximated and compares the two approaches (C2 and C3).

\section{Background}
\label{sec:background}

\subsection{Causal Block Diagrams}
\label{sec:cbd_background}
\ifnotreport
	\vspace{-5pt}
\fi

\begin{figure}[htb]
\begin{center}
  \includegraphics[width=0.8\textwidth]{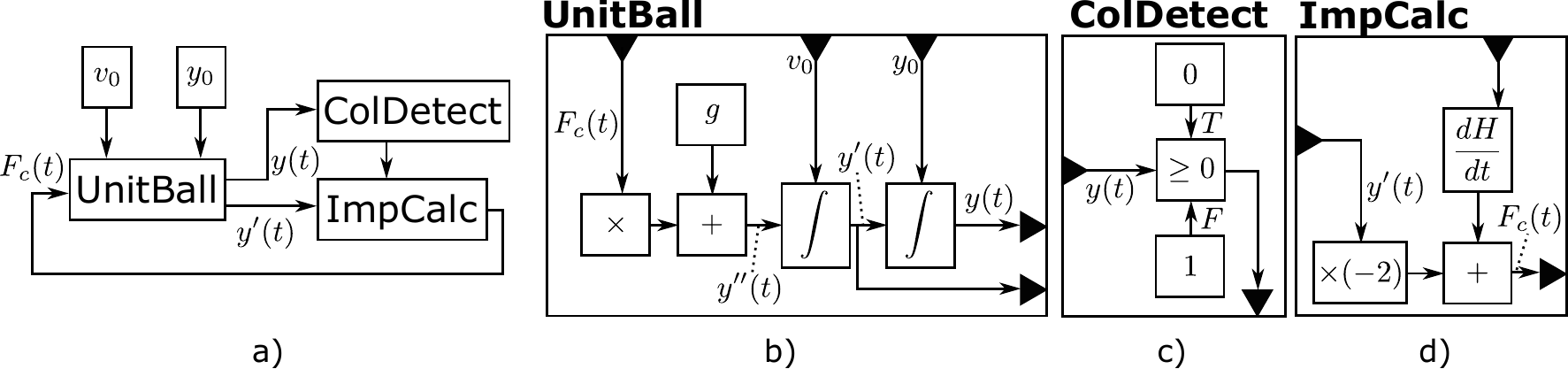}
  \caption{CBD model of a Bouncing Ball with unit mass (left), and specifications of each block (right).}
  \label{fig:bouncing_ball}
\end{center}
\end{figure}

Causal Block Diagrams (CBD) is a formalism used to model causal systems \shortcite{Posse2002,Gomes2016a}, commonly used in tools such as Simulink\textsuperscript{\textregistered}.
A CBD is expressed with three main entities: Blocks, Links and Ports.
Blocks denote operations or refer to other CBDs, defined elsewhere.
A Block may contain multiple input ports and a single output port, except if it represents a CBD, which then may contain any number of output ports.
Links establish the data-flow by connecting output ports to input ports.
\cref{fig:bouncing_ball} a) shows an example of a CBD that models a bouncing ball. 
The CBD on the left has 6 blocks, 3 representing the CBDs defined in \cref{fig:bouncing_ball} b), c) and d).

\ifnotreport
	\vspace{-5pt}
\fi

The basic procedure for the simulation of a CBD model is show in \algoref{alg:cbd_simulator} \shortcite{Mustafiz2016a,Vangheluwe2014}.
For the sake of simplicity, this procedure assumes that there are no algebraic loops in the model\footnote{Algebraic loops are supported by the provided simulator though.} \footnote{In \cref{fig:bouncing_ball} the integrator blocks only depend on their previous inputs, therefore no algebraic loop is formed.}.

\begin{algorithm}[htb]
	\begin{compactenum}
		\item Flatten CBD by replacing all blocks that refer to CBDs by their specification.
		\item Repeat until the end of the simulation:
			\begin{compactenum}
			\item Compute dependency graph between blocks and topologically sort the blocks using the dependency graph;
			\item For each block, compute the outputs from its inputs;
			\item Advance simulated time by $h$ units of time;
			\end{compactenum}
	\end{compactenum}
	\caption{CBD Simulator.}
	\label{alg:cbd_simulator}
\end{algorithm}
\ifnotreport
	\vspace{-5pt}
\fi

\subsection{Distributions}
\label{sec:distributions}
\ifnotreport
	\vspace{-5pt}
\fi

\begin{figure}[htb]
\centering
\begin{minipage}{.27\textwidth}
	\begin{center}
	  \includegraphics[width=1\textwidth]{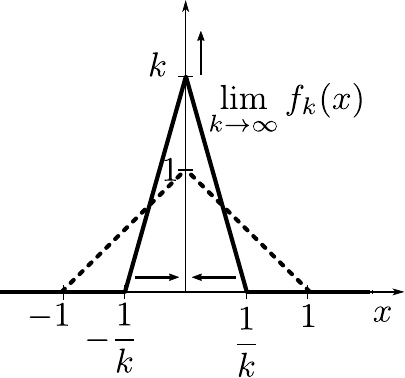}
	  \caption{Limit approximation of a Dirac delta.}
	  \label{fig:dirac_delta_approx}
	\end{center}
\end{minipage}%
\hspace{2pt}
\begin{minipage}{.66\textwidth}
	\begin{center}
	  \includegraphics[width=1\textwidth]{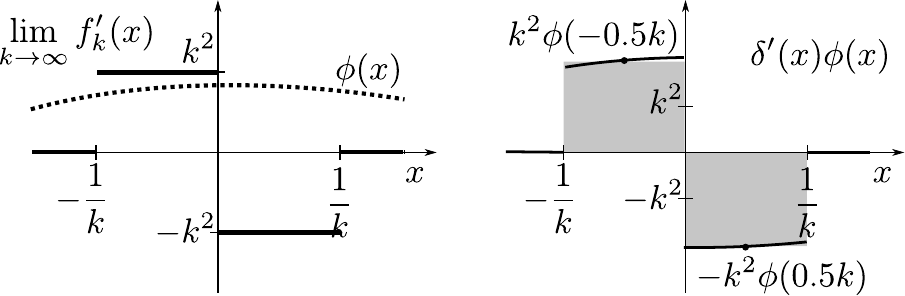}
	  \caption{Intuitive visualization of \cref{eq:delta_derivative_derivation}.}
	  \label{fig:dirac_delta_derivative}
	\end{center}
\end{minipage}
\end{figure}

Consider any function $f_k(x)$ that satisfies the conditions:
\begin{aligneq}\label{eq:dirac_properties}
	\pargroup{f_k(x) = 0 \text{ if } x \leq -\frac{1}{k} \text{ or } \frac{1}{k} \leq x}  
	\text{\ \ and } 
	\pargroup{\int_{-1/k}^{1/k} f_k(x) dx = 1}
\end{aligneq}
An example can be
$
f_k(x) = \conditional{\max(0, k + k^2 x)}{\max(0, k - k^2 x)}{x \leq 0}
$
plotted in dash stroke in \cref{fig:dirac_delta_approx} for $k=1$.

A Dirac delta $\delta(x)$ function can be constructed in a way that obeys to \cref{eq:dirac_properties} by taking the limit \shortcite{Strichartz2003}, depicted in \cref{fig:dirac_delta_approx}:
\begin{aligneq}\label{eq:delta_construction}
\delta(x) = \lim_{k \to \infty} f_k(x)
\end{aligneq}
One can informally say that
$$
\delta(x) = \conditional{0}{\infty}{x \neq 0}
$$
but this is as far as the utility of the statement goes.
A more interesting question is what happens when $\delta(x)$ is combined with a smooth (that is, infinitely differentiable) function $\phi(x)$ within an integral, that is, $\inftyint{\delta(x) \phi(x) dx}$.
According to Equations \ref{eq:dirac_properties} and \ref{eq:delta_construction}, and since $\phi$ is smooth, we have that $\delta(x) \phi(x) = \delta(x) \phi(0)$.
\ifreport
	To see why, let 
	$$
	\varphi_k(x)= \conditional{\phi(x)}{0}{-\frac{1}{k} < x < \frac{1}{k}}
	$$
	Then:
	$$
	\delta(x) \phi(x) = \lim_{k \to \infty} (f_k(x) \phi(x)) = \lim_{k \to \infty} (f_k(x) \varphi_k(x)) = \lim_{k \to \infty} f_k(x) \varphi_k(0)
	$$
\fi
The observed behavior then becomes:
\ifreport
	\begin{aligneq}\label{eq:dirac_distribution}
	\inftyint{\delta(x) \phi(x) dx} 
		& = \inftyint{\delta(x) \phi(0) dx}
		  = \phi(0) \inftyint{\delta(x) dx} \\
		&\ = \phi(0) \brackets{ \lim_{k \to \infty} \int_{-\frac{1}{k}}^{\frac{1}{k}} f_k(x) dx}
		  = \phi(0)
	\end{aligneq}
\else
	\begin{aligneq}\label{eq:dirac_distribution}
	\inftyint{\delta(x) \phi(x) dx} 
		& = \inftyint{\delta(x) \phi(0) dx}
		  = \phi(0) \inftyint{\delta(x) dx} 
		  = \phi(0)
	\end{aligneq}
\fi
which tells us that the Dirac delta $\delta(x-n)$ can be used to sample a signal $\phi$ at point $n$.

Consider now what happens when $\delta(x)$ is differentiated, depicted in \cref{fig:dirac_delta_derivative}:
\ifreport
	\begin{aligneq}\label{eq:delta_derivative_derivation}
	\inftyint{\dert{\delta}(x) \phi(x) dx} 
		&=  \lim_{k \to \infty} \int_{-1/k}^{1/k} \frac{f_k(x)-f_k(x-1/k)}{1/k}\phi(x) dx \\
		&=  \lim_{k \to \infty}
			\int_{-1/k}^0 \frac{f_k(0)-f_k(-1/k)}{1/k}\phi(x) dx 
			+ \int_0^{1/k} \frac{f_k(1/k)-f_k(0)}{1/k}\phi(x) dx  \\
		&=  \lim_{k \to \infty}
			\brackets{
				\int_{-1/k}^0 k^2 \phi(x) dx 
				+ \int_0^{1/k} -k^2 \phi(x) dx
			}  \\
		&=  \lim_{k \to \infty}
			\brackets{
				k^2 \phi(-\frac{1}{2k}) \frac{1}{k}
				- k^2 \phi(\frac{1}{2k}) \frac{1}{k}
			}  \\
		&=  \lim_{k \to \infty}
			\brackets{
				- \frac{\phi(\frac{1}{2k}) - \phi(-\frac{1}{2k})}{1/k}
			}  \\
		&= - \dert{\phi}(0)
	\end{aligneq}
\else
	\begin{aligneq}\label{eq:delta_derivative_derivation}
	\inftyint{\dert{\delta}(x) \phi(x) dx}
		&=  \lim_{k \to \infty}
			\brackets{
				k^2 \phi(-\frac{1}{2k}) \frac{1}{k}
				- k^2 \phi(\frac{1}{2k}) \frac{1}{k}
			}
		= - \dert{\phi}(0)
	\end{aligneq}
\fi
One can generalize the derivative of a Dirac delta to:
\begin{aligneq}\label{eq:general_dirac_derivative}
\inftyint{\dernt{\delta}{i}(x) \phi(x) dx}  = (-1)^i \dernt{\phi}{i}(0)
\end{aligneq}

These examples show that, even though $\delta(x)$ is not a function in the classical sense, one can observe interesting behavior, and identify an impulse from its interaction with other functions.
This is the essential idea of distributions.
A distribution is a function identified by the way it interacts with other \emph{test} functions \shortcite{Horvath1970,Friedlander1998,Strichartz2003}, and not by its direct plot (as with classical functions).

Any function $f(x)$ for which $\inftyint{|f(x)| dx} < \infty$, can be a distribution by writing it as $\tuple{f,\phi}$, with
\begin{aligneq}
\tuple{f,\phi} \triangleq \inftyint{ f(x) \phi(x) dx}
\text{ \ \ \ for any test function } \phi
\end{aligneq}
A test function is smooth and becomes 0 outside a bounded interval.

Two distributions are equal if \emph{the result of their interactions with any smooth function is the same}, that is,
$
f = g \iff \tuple{f,\phi} = \tuple{g,\phi} \text{ for all } \phi
$.
A consequence is that the two sides of \cref{eq:delta_construction} denote the same distribution.
This notion of equality allows badly behaved functions (discontinuous and/or with impulses) to exist, as long as their integral is finite.
As exemplified before, test functions play a key role  because they tend to compensate for the bad shape of the distribution (recall \cref{eq:general_dirac_derivative}). 

Consider the Heaviside distribution $H$ defined by:
\begin{aligneq}\label{eq:heaviside_distribution}
H(x) = \conditional{1}{0}{x \geq 0}
\end{aligneq}
Even though it has a discontinuity at the origin, it can be differentiated as a distribution:
\begin{aligneq}
\tuple{\dert{H},\phi} 
	= \inftyint{ \dert{H}(x) \phi(x) dx}
	= \brackets{H(x)\phi(x)}_{-\infty}^\infty - \inftyint{ H(x) \dert{\phi}(x) dx} 
\ifreport
	\\
\fi
	= - \int_0^\infty \dert{\phi}(x) dx
	= - \brackets{\phi(x)}_{0}^\infty 
	= \phi(0)
\end{aligneq}
According to \cref{eq:dirac_distribution}, $\tuple{\delta,\phi}=\phi(0)$. 
So, by distribution equality,
\begin{aligneq}\label{eq:heaviside_dirac_relation}
\dert{H} = \delta
\end{aligneq}

To see the relationship between distributions and impulsive differential equations, recall the bouncing ball model in \cref{eq:bouncing_ball_dirac}, and suppose that its solution is given by
$$
y(t) = \underbrace{\pargroup{y_0 + v_0 t - \frac{1}{2} g t^2}}_{\text{before collision}}(1 - H(t-t_c)) 
	+ \underbrace{\pargroup{y(t_c^-) - \dert{y}(t_c^-)(t-t_c)  - \frac{1}{2} g (t-t_c)^2}}_{\text{after collision}}H(t-t_c)
$$
where $t_c$ is the time at which the collision occurs.
\ifreport
	This solution is the mix of two trajectories: the one before the bounce, and the one after.
\fi

It is clear that derivatives of the solution are discontinuous functions, so differentiating $y(t)$ twice and considering the result as a distribution yields:
\ifreport
	\begin{aligneq}
	\tuple{\dernt{y}{2},\phi} 
		&= \inftyint{\dernt{y}{2}(t)\phi(t) dt} 
		= \inftyint{y(t)\dernt{\phi}{2}(t) dt} \\
		&= \inftyint{
			\pargroup{y_0 + v_0 t - \frac{1}{2} g t^2}
			(1 - H(t-t_c))
			\dernt{\phi}{2}(t) dt
		}\\
		& \hspace{3em}
		+ \inftyint{
			 \pargroup{y(t_c^-) - \dert{y}(t_c^-)(t-t_c)  - \frac{1}{2} g (t-t_c)^2}
			 H(t-t_c) 
			\dernt{\phi}{2}(t) dt
		} \\
		&= \int_{-\infty}^{t_c} \pargroup{y_0 + v_0 t - \frac{1}{2} g t^2} \dernt{\phi}{2}(t) dt \\
		& \hspace{3em}
		+ \int_{t_c}^\infty \pargroup{y(t_c^-) - \dert{y}(t_c^-)(t-t_c)  - \frac{1}{2} g (t-t_c)^2} \dernt{\phi}{2}(t) dt   \\
		&= \brackets{\pargroup{y_0 + v_0 t - \frac{1}{2} g t^2} \dert{\phi}(t)}_{-\infty}^{t_c} 
		- \pargroup{ \brackets{\pargroup{v_0 t - g t} \phi(t)}_{-\infty}^{t_c} - \int_{-\infty}^{t_c} - g \phi(t) dt}
		\\
		& \hspace{3em}
		+ \brackets{\pargroup{y(t_c^-) -  \dert{y}(t_c^-)(t-t_c)  - \frac{1}{2} g (t-t_c)^2} \dert{\phi}(t)}_{t_c}^{\infty}
		\\
		& \hspace{3em}
		- \pargroup{
		\brackets{\pargroup{- \dert{y}(t_c^-) - g(t-t_c)} \phi(t)}_{t_c}^{\infty}
		- \int_{t_c}^\infty - g \phi(t) dt
		} 
		\\
		&= y(t_c^-)\dert{\phi}(t_c)
			-  \dert{y}(t_c^-) \phi(t_c) 
			- g \int_{-\infty}^{t_c} \phi(t) dt \\
			& \hspace{3em}
			- y(t_c^-)\dert{\phi}(t)
			- \dert{y}(t_c^-)\phi(t_c)
			- g \int_{t_c}^\infty  \phi(t) dt \\
		&= - 2 \dert{y}(t_c^-) \phi(t_c) 
			- g \int_{-\infty}^{\infty} \phi(t) dt \\
		&= - 2 \dert{y}(t_c^-) \int_{-\infty}^{\infty} \delta(t-t_c) \phi(t) dt 
			- g \int_{-\infty}^{\infty} \phi(t) dt \\
		&= 	\int_{-\infty}^{\infty} \pargroup{-g - 2 \dert{y}(t_c^-) \delta(t-t_c)} \phi(t) dt 
			= \tuple{-g - 2 \dert{y}(t_c^-) \delta(t-t_c), \phi(t)}
	\end{aligneq}
\else
	\begin{aligneq}
	& \tuple{\dernt{y}{2},\phi}
		= \inftyint{\dernt{y}{2}(t)\phi(t) dt} 
		= \inftyint{y(t)\dernt{\phi}{2}(t) dt} \\
	&\ \ = \int_{-\infty}^{t_c} \pargroup{y_0 + v_0 t - \frac{1}{2} g t^2} \dernt{\phi}{2}(t) dt
		+ \int_{t_c}^\infty \pargroup{y(t_c^-) - \dert{y}(t_c^-)(t-t_c)  - \frac{1}{2} g (t-t_c)^2} \dernt{\phi}{2}(t) dt   \\
	&\ \ = \tuple{-g - 2 \dert{y}(t_c^-) \delta(t-t_c), \phi(t)}
	\end{aligneq}
\fi
which, according ot the distribution equality, obeys to \cref{eq:bouncing_ball_dirac}.
In the next section we apply this theory for the derivation of the symbolic computations on impulses.

\ifnotreport
	\vspace{-10pt}
\fi
\section{Symbolic Manipulation of Dirac Deltas}
\label{sec:hybridcbds:semantics}

In the formalization of Dirac CBDs, the signals transmitted in the links between blocks are distributions that can have discontinuities and impulses (and any derivative of the impulses).
Formally, let $\set{\tau_j} \subset \setreal$ denote the sequence of all times at which an impulse (or any derivative of it) occurs. 
Then, following the formalization in \cite{Nilsson2003}, the signals are of the form:
\begin{aligneq}\label{eq:generalized_signal_representation}
S(t) = s(t) + \sum_{i=0}^n \sumseq{\tau_j} a_{ij} \dernt{\delta}{i}( t-\tau_j) 
\end{aligneq}
where $s$ denotes a \emph{piece-wise continuous} impulse-free function, $n$ is the maximum derivative order of any impulse in the signal and $a_{ij}$ is a (possibly zero) constant called the impulse coefficient.
At any time $\tau_j$ impulses occur, both the right and left limits of $S(t)$ have to have the same impulses. 
This ensures that the derivations in this section can disregard the left and right limits of the impulse part of the signal. 
Formally, 
\begin{aligneq}\label{eq:generalized_signal_conditions}
S(t^+) - S(t^-) = s(t^+) - s(t^-)
\end{aligneq}
Discontinuities in the impulse free part of the signal, that is, $s(t^+) \neq s(t^-)$, are allowed.

\begin{figure}[htb]
\begin{center}
	\includegraphics[width=0.5\textwidth]{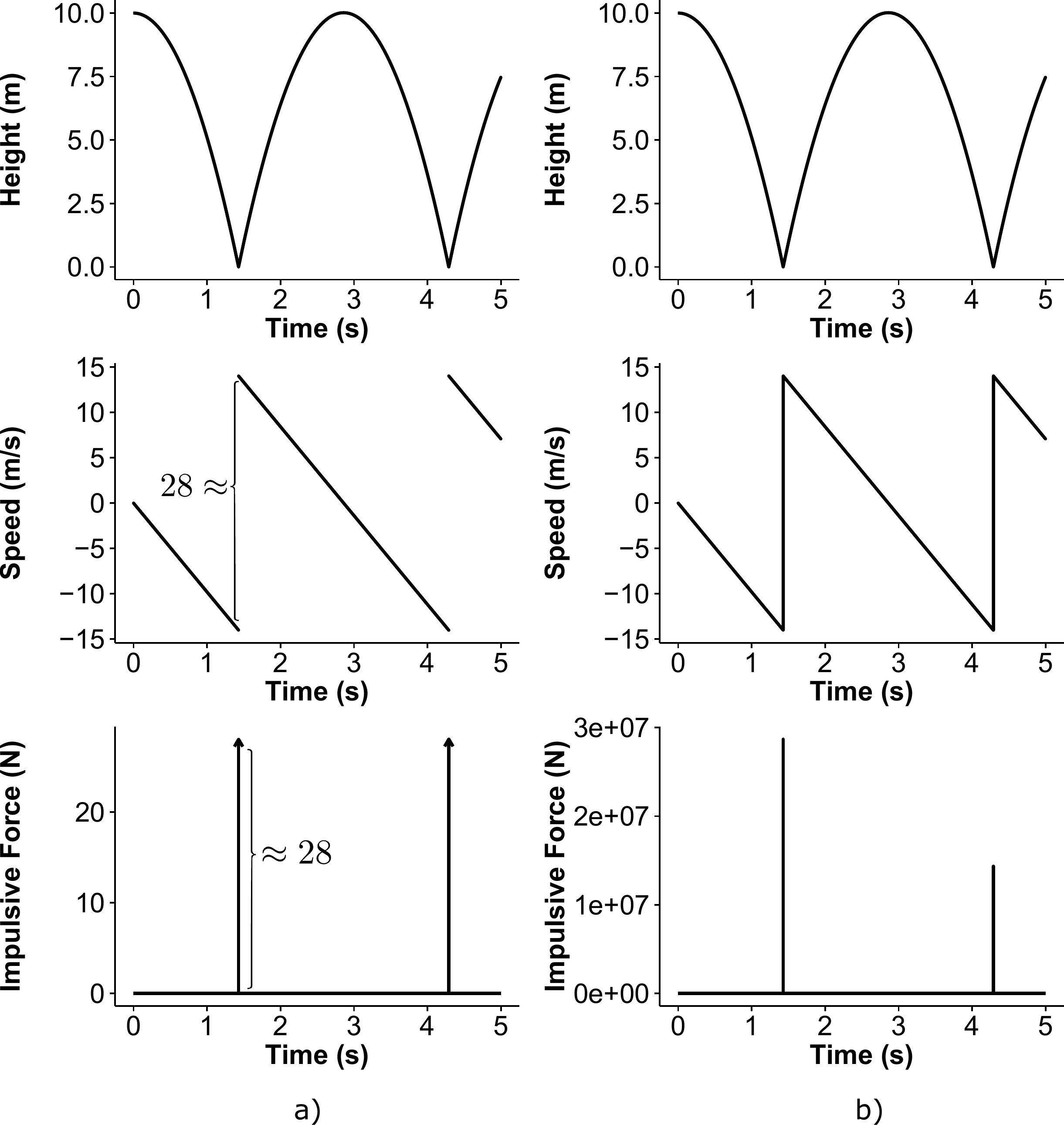}
	\caption{Bouncing ball trace, with perfectly elastic collision. The symbolic-impulse approach is a) and the approximated-impulse one is b).}
	\label{fig:sample_bouncing_ball_trace}
\end{center}
\end{figure}

The CBD blocks denote manipulations on distributions in the form of \cref{eq:generalized_signal_representation}.
To ensure that these are compositional, we need to show that the output $Y(t)$ of each block $b$ is of the form of  \cref{eq:generalized_signal_representation}, assuming that the inputs are too.
Furthermore, to conform strictly to the theory, the moment any signal is plotted, it should be plotted as $\tuple{S(t), \varphi(t)}$, for suitable test function.
However, the utility of this simulation package is that it allows the modeller to observe the discontinuities and impulses in a signal, so any impulse is plotted with an arrow. The height of the arrow represents the impulse coefficient.
An example is shown in \cref{fig:sample_bouncing_ball_trace} a).
We now provide a derivation of each operation.
\ifnotreport
	For more details, see \shortcite{Gomes2016c}.
\fi

The \textbf{Sum} block takes two signals ($U(t)$ and $V(t)$) as inputs and produces $Y(t)$ as output:
\ifreport
	\begin{aligneq}\label{eq:adder_output_derivation}
	Y(t) &= U(t) + V(t) \Leftrightarrow
	\tuple{Y(t),\varphi(t)} = \tuple{U(t) + V(t),\varphi(t)} \text{ \ \ \ for any test function } \phi \\
				&=\tuple{U(t),\varphi(t)} + \tuple{V(t),\varphi(t)} \text{ by linearity}\\	
				&=\tuple{u(t) + \sum_{i=0}^{n_u} \sumseq{\tau^u_j} a_{ij} \dernt{\delta}{i}(t-\tau^u_j),\varphi(t)} + \tuple{v(t) + \sum_{i=0}^{n_v} \sumseq{\tau^v_j} b_{ij} \delta^{(i)}( t-\tau^v_j),\varphi(t)}\\
				&=\tuple{u(t),\varphi(t)} + \tuple{v(t),\varphi(t)} + \sum_{i=0}^{n_u} \sumseq{\tau^u_j} a_{ij} \tuple{\dernt{\delta}{i}(t-\tau^u_j),\varphi(t)} + \\
				& \hspace{2em} \sum_{i=0}^{n_v} \sumseq{\tau^v_j} b_{ij} \tuple{\delta^{(i)}( t-\tau^v_j),\varphi(t)}  \\
				& \text{ by inverse linearity} \\
				&=\tuple{u(t) + v(t) + \sum_{i=0}^{n_u} \sumseq{\tau^u_j} a_{ij} \dernt{\delta}{i}(t-\tau^u_j) + \sum_{i=0}^{n_v} \sumseq{\tau^v_j} b_{ij} \delta^{(i)}( t-\tau^v_j),\varphi(t)}
	\end{aligneq}
\else
	\begin{aligneq}\label{eq:adder_output_derivation}
	& Y(t) = U(t) + V(t) \Leftrightarrow
	\tuple{Y(t),\varphi(t)} = \tuple{U(t) + V(t),\varphi(t)} \text{ \ \ \ for any test function } \phi \\
	& \ \ =\tuple{U(t),\varphi(t)} + \tuple{V(t),\varphi(t)} \\
	& \ \ =\tuple{u(t) 
					+ v(t) 
					+ \sum_{i=0}^{n_u} \sumseq{\tau^u_j} a_{ij} \dernt{\delta}{i}(t-\tau^u_j) 
					+ \sum_{i=0}^{n_v} \sumseq{\tau^v_j} b_{ij} \delta^{(i)}( t-\tau^v_j),\varphi(t)}
	\end{aligneq}
\fi
Operationally, it adds the two impulse-free signals $u(t)$ and $v(t)$, and adds the coefficients $a_{ij}$ of the impulses in $U(t)$, with the coefficients $b_{ij}$ of the impulses in $V(t)$.
$Y(t)$ is of the form of \cref{eq:generalized_signal_representation} so the derivation is complete.

The \textbf{Negation} block is derived in the same way. It negates the impulse-free part of the signal and all the coefficients.
\ifreport
	\begin{aligneq}
	\tuple{Y(t),\varphi(t)} &= \tuple{-U(t),\varphi(t)} \\
				&= \tuple{-\pargroup{u(t) + \sum_{i=0}^{n_u} \sumseq{\tau_j} a_{ij} \dernt{\delta}{i}(t-\tau_j)},\varphi(t)}\\
				&= -\tuple{u(t),\varphi(t)} - \sum_{i=0}^{n_u} \sumseq{\tau_j} a_{ij} \tuple{\dernt{\delta}{i}(t-\tau_j),\varphi(t)}\\
				&= \tuple{-u(t) - \sum_{i=0}^{n_u} \sumseq{\tau_j} a_{ij} \dernt{\delta}{i}(t-\tau_j) ,\varphi(t)}
	\end{aligneq}
\fi

The \textbf{Switch} block outputs either $0$ or $1$, according to a condition on the input $C(t)$. 
A family of pathological cases related to impulses can be identified if $C(t)$ is allowed to have impulses. 
For example, does a function which at a point $t_0>0$ is positive, but has a negative impulse (e.g., $C(t_0) = c(t_0)^2 - \delta(0)$) cross the zero?
So $C(t)$ is assumed to be free of impulses, that is, $C(t) = c(t)$. 
The derivation is then:
\begin{aligneq}\label{eq:switch_block_io}
\tuple{Y(t),\varphi(t)} = \tuple{H(c(t)),\varphi(t)} 
\end{aligneq}
where $H(x)$ is the Heaviside distribution defined in \cref{eq:heaviside_distribution}.
\ifreport
	Alternatively, it can be written as:
	\begin{aligneq}
	\tuple{Y(t),\varphi(t)} = \tuple{\conditional{1}{0}{c(t) \geq 0} ,\varphi(t)} 
	\end{aligneq}
	Note however, that $c(t)$ can still be piece-wise continuous.
	Then the left limit of the output, denoted as $Y(t^-)$, is evaluated with $c(t^-) > 0$ and the right limit $Y(t^+)$ with $c(t^+) \geq 0$. 
	In other words, for any time $t_0$ where $c(t_0^-) < 0$ and $c(t_0^+) > 0$, $H(c(t_0^-)) = 0$ and $H(c(t_0^+)) = 1$.
\fi

The \textbf{Decision} block is a generalization of the Switch block. It forwards one of two inputs ($U(t)$ and $V(t)$) to the output $Y(t)$ according to an input condition $C(t) \geq 0$.
Similarly to the Switch block, we have to assume that $C(t) = c(t)$ but in addition, at any time $t_0$ where $c(t_0^-) < 0$ and $c(t_0^+) \geq 0$, there can be no impulses on either $U(t_0)$ or $V(t_0)$. 
If this were allowed, a signal could be created which violates \cref{eq:generalized_signal_conditions}.
Under these assumptions, the derivation is:
\begin{aligneq}
\tuple{Y(t),\varphi(t)} = \tuple{U(t)H(C(t)) + V(t)(1 - H(C(t))),\varphi(t)} 
\end{aligneq}

\ifreport
	The output of the Decision block can also be written as
	\begin{aligneq}
	\tuple{Y(t),\varphi(t)} = \tuple{\conditional{U(t)}{V(t)}{c(t) \geq 0} ,\varphi(t)} 
	\end{aligneq}
\fi

The \textbf{Derivative} block is partially derived as follows:
\ifreport
	\begin{aligneq}\label{eq:derivative_functional_partial}
	\tuple{Y(t),\varphi(t)} &= \tuple{\dert{U}(t),\varphi(t)} \\
		&= \tuple{\dert{u}(t) + \dert{\brackets{\sum_{i=0}^{n_u} \sumseq{\tau_j} a_{ij} \dernt{\delta}{i}(t-\tau_j)(t)}},\varphi(t)} \\
		&= \tuple{\dert{u}(t),\varphi(t)} + \tuple{\sum_{i=0}^{n_u} \sumseq{\tau_j} a_{ij} \dernt{\delta}{i+1}(t-\tau_j)(t),\varphi(t)}
	\end{aligneq}
\else
	\begin{aligneq}\label{eq:derivative_functional_partial}
	\tuple{Y(t),\varphi(t)} &= \tuple{\dert{U}(t),\varphi(t)}
		= \tuple{\dert{u}(t),\varphi(t)} + \tuple{\sum_{i=0}^{n_u} \sumseq{\tau_j} a_{ij} \dernt{\delta}{i+1}(t-\tau_j)(t),\varphi(t)}
	\end{aligneq}
\fi
The term $\tuple{\dert{u}(t),\varphi(t)}$ must be further developed because it can have discontinuities.
\ifreport
For example, when the Derivative block is connected to the output of a Decision block.
\fi 
Let $\set{t_d}$ be the countable sequence of times at which $u(t_d^-) \neq u(t_d^+)$.
In the proximity of a discontinuity at time $t_d$, $u(t)$ is described as 
$u(t) = u(t^-)(1-H(t-t_d)) + u(t^+)H(t-t_d)$. Intuitively, it can be seen as the output signal of a Decision block.
$\tuple{\dert{u}(t),\varphi(t)}$ is then:
\ifreport
	\begin{aligneq}
	\tuple{\dert{u}(t),\varphi(t)} &= 
			\tuple{\dertb{
				u(t^-)(1-H(t-t_d)) + u(t^+)H(t-t_d)
				},\varphi(t)} \\
		&= \tuple{
			\dertb{u(t^-) 
					+ (u(t^+) - u(t^-))H(t-t_d)
				}
			,\varphi(t)} \\
		&= \tuple{
			\dert{u}(t^-)
				+ \dertb{(u(t^+) - u(t^-))H(t-t_d)}
			,\varphi(t)} \\
		&= \tuple{\dert{u}(t^-),\varphi(t)} 
			+ \tuple{\dertb{(u(t^+) - u(t^-))H(t-t_d)},\varphi(t)}
	\end{aligneq}
	The distribution $\tuple{\dertb{(u(t^+) - u(t^-))H(t-t_d)},\varphi(t)}$ can be further simplified:
	\begin{aligneq}
	&\tuple{\dertb{(u(t^+) - u(t^-))H(t-t_d)},\varphi(t)} =
	 - \tuple{(u(t^+) - u(t^-))H(t-t_d),\dert{\varphi}(t)} \\
		&\hspace{3em}=  - \inftyint{(u(t^+) - u(t^-))H(t-t_d)\dert{\varphi}(t)} dt \\
		&\hspace{3em}=  - \int_{t_d}^{\infty}(u(t^+) - u(t^-))\dert{\varphi}(t) dt \\
		&\hspace{3em} \text{(integration by parts)} \\
		&\hspace{3em}=  - \pargroup{ \brackets{(u(t^+) - u(t^-))\varphi(t)}_{t_d}^\infty} - \int_{t_d}^{\infty}(\dert{u}(t^+) - \dert{u}(t^-))\varphi(t) dt \\
		&\hspace{3em}=   (u(t_d^+) - u(t_d^-))\varphi(t_d)  + \int_{-\infty}^{\infty}(\dert{u}(t^+) - \dert{u}(t^-))H(t-t_d)\varphi(t) dt \\
		&\hspace{3em}=   (u(t_d^+) - u(t_d^-))
				\tuple{\delta(t-t_d), \varphi(t)} 
			+ \tuple{(\dert{u}(t^+) - \dert{u}(t^-))H(t-t_d), \varphi(t)} \\
	\end{aligneq}
	In conclusion, around each discontinuity $t_d \in T_d$, we have:
	{\small
	\begin{aligneq}\label{eq:derivative_discontinuity}
	\tuple{\dert{u}(t),\varphi(t)} &= 	
		\tuple{
			\dert{u}(t^-)(1 - H(t-t_d))
			+ \dert{u}(t^+)H(t-t_d) 
			+ (u(t_d^+) - u(t_d^-))\delta(t-t_d)
		,\varphi(t)}
	\end{aligneq}
	}
\else
	\begin{aligneq}\label{eq:derivative_discontinuity}
	&\tuple{\dert{u}(t),\varphi(t)} 
	= \tuple{\dert{u}(t^-),\varphi(t)} 
			+ \tuple{\dertb{(u(t^+) - u(t^-))H(t-t_d)},\varphi(t)} \\
	&\ \ = \tuple{\dert{u}(t^-),\varphi(t)} 
		- \tuple{(u(t^+) - u(t^-))H(t-t_d),\dert{\varphi}(t)} 
	= \tuple{\dert{u}(t^-),\varphi(t)} 
		- \int_{t_d}^{\infty}(u(t^+) - u(t^-))\dert{\varphi}(t) dt  \\
	&\ \ = \tuple{
		\dert{u}(t^-)(1 - H(t-t_d))
		+ \dert{u}(t^+)H(t-t_d) 
		+ (u(t_d^+) - u(t_d^-))\delta(t-t_d)
	,\varphi(t)}
	\end{aligneq}
\fi
This is interpreted as originating an impulse, with a coefficient given by the magnitude of the discontinuity. 
The derivative signal can itself be discontinuous, if its left limit is also different than its right limit. 
Across all possible discontinuities $t_d \in T_d$ of $u(t)$, and abstracting the discontinuities in  $\dert{u}(t)$, the signal in \cref{eq:derivative_discontinuity} is in the form of \cref{eq:generalized_signal_representation}.
\ifreport
	Concretely:
	\begin{aligneq}\label{eq:derivative_discontinuity_full}
	\tuple{\dert{u}(t),\varphi(t)} &= 	
		\tuple{
			\dert{u}(t) 
			+ \sum_{t_d \in T_d} (u(t_d^+) - u(t_d^-))\delta(t-t_d)
		,\varphi(t)}
	\end{aligneq}
\fi
Joining \cref{eq:derivative_discontinuity} and \cref{eq:derivative_functional_partial} yields:
\ifreport
	\small
\fi
\begin{aligneq}
\tuple{Y(t),\varphi(t)} &= 
	\tuple{
	\dert{u}(t) 
		+ \sum_{t_d \in \set{t_d}} (u(t_d^+) - u(t_d^-))\delta(t-t_d)
		+ \sum_{i=0}^{n_u} \sumseq{\tau_j} a_{ij} \dernt{\delta}{i+1}(t-\tau_j)(t)
	,\varphi(t)
	}
\end{aligneq}
\ifreport
	\normalsize
\fi

\ifreport
	For example, suppose the input to the Derivative block is the signal\\
	$$
	U(t)= (v_0 - gt)(1-H(t-t_d)) + (- (v_0 - gt_d) - g(t-t_d))H(t-t_d)
	$$
	With $v_0$,$g$ constants.
	Then the output of the Derivative block is
	$$\tuple{Y(t),\varphi(t)} = -g  -2 (v_0 - gt_d)\delta(t-t_d)$$ 
\fi

The \textbf{Integrator} block is derived as:
\ifreport
	\begin{aligneq}
	\tuple{Y(t),\varphi(t)} &= \tuple{\int_0^t U(x) dx,\varphi(t)} \\
				&= \tuple{\int_0^t u(x) + \sum_{i=0}^{n_u} \sumseq{\tau_j} a_{ij} \dernt{\delta}{i}(x-\tau_j) dx,\varphi(t)} \\
				&= \tuple{\int_0^t u(x) dx, \varphi(t)} + 
					\sum_{i=0}^{n_u} \sumseq{\tau_j} a_{ij} \tuple{\int_0^t \dernt{\delta}{i}(x-\tau_j) dx, \varphi(t)} \\
				&= \tuple{\int_0^t u(x) dx, \varphi(t)} 
				+ \sumseq{\tau_j} a_{0j} \tuple{\int_0^t \delta(x-\tau_j) dx, \varphi(t)} 
				\\
				&\hspace{3em}+ \sum_{i=1}^{n_u} \sumseq{\tau_j} a_{ij} \tuple{\dernt{\delta}{i-1}(t-\tau_j), \varphi(t)} \\
				&= \tuple{\int_0^t u(x) dx, \varphi(t)} 
				+ \sumseq{\tau_j} a_{0j} \tuple{H(x-\tau_j), \varphi(t)} 
				\\
				&\hspace{3em}+ \sum_{i=1}^{n_u} \sumseq{\tau_j} a_{ij} \tuple{\dernt{\delta}{i-1}(t-\tau_j), \varphi(t)} \\
				&= \tuple{
					\int_0^t u(x) dx
					+ \sumseq{\tau_j} a_{0j} H(x-\tau_j)
					+ \sum_{i=1}^{n_u} \sumseq{\tau_j} a_{ij} \dernt{\delta}{i-1}(t-\tau_j) 
					, \varphi(t)} \\			
	\end{aligneq}
\else
	\begin{aligneq}
	& \tuple{Y(t),\varphi(t)} = \tuple{\int_0^t U(x) dx,\varphi(t)}
		= \tuple{
					\int_0^t u(x) dx
					+ \sumseq{\tau_j} a_{0j} H(x-\tau_j)
					+ \sum_{i=1}^{n_u} \sumseq{\tau_j} a_{ij} \dernt{\delta}{i-1}(t-\tau_j) 
					, \varphi(t)} \\			
	\end{aligneq}
\fi
Which matches the intuition: the Integrator block computes the normal integration of the impulse-free signal, reduces the order of any impulse derivative, and computes a discontinuity whenever an impulse is integrated (recall \cref{eq:heaviside_dirac_relation}).
The magnitude of the discontinuity is the impulse coefficient.

The \textbf{Product} of two distributions is not defined in general. 
However, if one of the input signals to the Product block is impulse-free\footnote{This restriction cannot be enforced automatically by the tool.}, then it can be derived as:
\ifreport
	\small
		\begin{aligneq}\label{eq:product_derivation}
		\tuple{Y(t),\varphi(t)} &= \tuple{U(t)V(t),\varphi(t)} \\
			&=\tuple{u(t) \pargroup{v(t) + \sum_{i=0}^{n_v} \sumseq{\tau_j} a_{ij} \dernt{\delta}{i}( t-\tau_j)},\varphi(t)} \text{ by assumption}\\	
			&= \tuple{u(t) v(t),\varphi(t)} + \tuple{\sum_{i=0}^{n_v} \sumseq{\tau_j} a_{ij} u(t) \dernt{\delta}{i}( t-\tau_j),\varphi(t)}\\	
			&= \tuple{u(t) v(t),\varphi(t)} + \sum_{i=0}^{n_v} \sumseq{\tau_j} a_{ij} \tuple{ u(t)  \dernt{\delta}{i}( t-\tau_j),\varphi(t)}\\	
			&  \text{ by derivative and product properties} \\
			&= \tuple{u(t)  v(t),\varphi(t)} + \sum_{i=0}^{n_v} \sumseq{\tau_j} a_{ij}  (-1)^i \tuple{ \delta( t-\tau_j), \dernt{\brackets{u(t)  \varphi(t)}}{i} } \\
			& \text{ by general Leibniz rule} \\
			&= \tuple{u(t)  v(t),\varphi(t)} + \sum_{i=0}^{n_v} \sumseq{\tau_j} a_{ij}  (-1)^i \tuple{ \delta( t-\tau_j), \sum_{k=0}^i \binom{i}{k} \dernt{u}{k}(t) \dernt{\varphi}{i-k}(t) }  \\
			&  \text{ by linearity and delta properties}  \\
			&= \tuple{u(t)  v(t),\varphi(t)} + \sum_{i=0}^{n_v} \sumseq{\tau_j} a_{ij}  (-1)^i \sum_{k=0}^i \binom{i}{k}  \dernt{u}{k}(\tau_j) \dernt{\varphi}{i-k}(\tau_j) \\
			&  \text{ by reverse delta property and } \dernt{u}{k}(\tau_j) \text{ constant}\\
			&= \tuple{u(t)  v(t),\varphi(t)} + \sum_{i=0}^{n_v} \sumseq{\tau_j} a_{ij}  (-1)^i \sum_{k=0}^i \binom{i}{k}  \dernt{u}{k}(\tau_j) \tuple{\delta( t-\tau_j), \dernt{\varphi}{i-k}(t)} \\
			&  \text{ by property } (-1)^{i-k}  (-1)^{i-k} = 1\\
			&= \tuple{u(t)  v(t),\varphi(t)} + \\
			&\hspace{3em}	\sum_{i=0}^{n_v} \sumseq{\tau_j} a_{ij}  (-1)^i \sum_{k=0}^i \binom{i}{k}  \dernt{u}{k}(\tau_j) (-1)^{i-k}  (-1)^{i-k} \tuple{\delta( t-\tau_j), \dernt{\varphi}{i-k}(t)} \\
			&  \text{ by reverse derivative property}\\
			&= \tuple{u(t)  v(t),\varphi(t)} + \sum_{i=0}^{n_v} \sumseq{\tau_j} a_{ij}  \sum_{k=0}^i \binom{i}{k}  \dernt{u}{k}(\tau_j)  (-1)^i  (-1)^{i-k} \tuple{\dernt{\delta}{i-k}( t-\tau_j), \varphi(t)} \\
			&  \text{ by fact } (-1)^i  (-1)^{i-k} = (-1)^k \\
			&= \tuple{u(t)  v(t),\varphi(t)} + \sum_{i=0}^{n_v} \sumseq{\tau_j} a_{ij}  \sum_{k=0}^i \binom{i}{k}  \dernt{u}{k}(\tau_j)  (-1)^k \tuple{\dernt{\delta}{i-k}( t-\tau_j), \varphi(t)}
		\end{aligneq}
	\normalsize
\else
	\begin{aligneq}\label{eq:product_derivation}
	\tuple{Y(t),\varphi(t)} 
	&= \tuple{U(t)V(t),\varphi(t)} \\
	&= \tuple{u(t)  v(t),\varphi(t)} + \sum_{i=0}^{n_v} \sumseq{\tau_j} a_{ij}  \sum_{k=0}^i \binom{i}{k}  \dernt{u}{k}(\tau_j)  (-1)^k \tuple{\dernt{\delta}{i-k}( t-\tau_j), \varphi(t)}
	\end{aligneq}
\fi

Essentially, each product between a function and a Dirac delta derivative yields a set of simultaneous impulses of decreasing derivative orders, that are multiplied with increasing derivatives of the function. 
The computation of the coefficients, therefore, requires that the derivatives of $u(t)$, order up to (and including) $n_v$, be known or estimated.
To get the intuition, recall \cref{eq:delta_derivative_derivation} and \cref{fig:dirac_delta_derivative}, and note that 
\begin{aligneq}\label{eq:delta_derivative_simple_derivation}
& \tuple{\dert{\delta}(t-\tau_j)u(t),\varphi(t)} 
  = -\dert{\brackets{u(t-\tau_j)\varphi(t-\tau_j)}}
  = \tuple{ u(t) \dert{\delta} - \dert{u}(t)\delta,\varphi(t)}
\end{aligneq}

\ifreport
	For example, if $U(t)=-g t$ and $V(t) = 1 + 20\dernt{\delta}{2}(t-1.44)$, then the Product block computes:
	\begin{aligneq}
	\tuple{U(t)V(t),\varphi(t)} &=  
		\tuple{- g t \times 1,\varphi(t)} + 
		20  \sum_{k=0}^2 \binom{2}{k}  \dernt{u}{k}(1.44)  (-1)^k \tuple{\dernt{\delta}{2-k}( t-1.44), \varphi(t)} \\ 
		&=  \tuple{- g t,\varphi(t)} - 28.8 g \tuple{\dernt{\delta}{2}( t-1.44), \varphi(t)}  +
			40 g \tuple{\dernt{\delta}{1}( t-1.44), \varphi(t)}
	\end{aligneq}
\fi

The \textbf{Inverter} block, to be well defined, requires that the input is free of impulses:
\begin{aligneq}
\tuple{Y(t),\varphi(t)} = \tuple{\frac{1}{u(t)},\varphi(t)}
\end{aligneq}

This section has provided the relationship between DiracCBDs and Distributions, and derived each block used in the CBD formalism.
In the DiracCBDs simulator the impulses are treated symbolically but the rest of the signals are discretized approximations of the continuous system.
In particular, if $h>0$ is the simulation time step (recall \algoref{alg:cbd_simulator}), then the integrator block can be approximated by recursive right Riemann's sum and the derivative block by finite difference: 
\begin{aligneq}\label{eq:numerical_approximations}
	Y_{\text{int}}(t) \approx Y_{\text{int}}(t-h) + U(t) \times h
	\text{\ \ and \ \ }
	Y_{\text{der}}(t) \approx \frac{U(t) - U(t-h)}{h}
\end{aligneq}
The signals, in the form of \cref{eq:generalized_signal_representation}, can be encoded in many ways.
For example, at each time step, the left and right limit of the signals can be stored so that discontinuities can be easily identified, and a matrix of impulse coefficients kept. 
Each block then takes the input left/right limit and the matrix, and produces an output of the same form.
\shortcite{Nilsson2003} encodes the left limit of the signal, along with all its derivatives (lazily evaluated).
Another example encoding can be seen in \shortcite{Lee2015} that is based in the super dense time encoding.
\cref{fig:sample_bouncing_ball_trace} a) shows the bouncing ball trace produced by our simulator.

\section{Numerical Approximation of Dirac Deltas}
\label{sec:general_cbds_discussion}

The previous section offered insight on how to symbolically manipulate impulses in the context of simulation. 
This manipulation can get complex and counter intuitive, especially with derivatives of impulses (recall \cref{eq:product_derivation} and \cref{eq:delta_derivative_simple_derivation}), and requires a special encoding of the signals.
The purpose of this section is to describe an alternative approach by approximating Dirac deltas as large functions, and see how it compares to the symbolic approach.

A Dirac delta impulse is a arbitrary high function at a very small interval of time that obeys to \cref{eq:dirac_properties}.
In a numerical simulation, the smallest interval is $h$, so a good numerical approximation of an impulse $\delta(t-\tau_d)$ would be to produce a large value $N$ at a simulation step that is immediately after $\tau_d$:
\begin{aligneq}
	\delta(t-\tau_d) \approx \conditional{N}{0}{0 \leq t - \tau_d < h}
\end{aligneq}
We cannot just pick an arbitrary high value for $N$ (e.g., $N = 2^{31} - 1$) because of the conditions in \cref{eq:dirac_properties}.

\cref{eq:heaviside_dirac_relation} holds the key to numerically approximate the impulse.
Suppose the Dirac delta is constructed as:
\begin{aligneq}\label{eq:alternative_construct_delta}
	\delta(x) = \lim_{k\to\infty} \dert{H_k}(x) = \lim_{k\to\infty} \conditional{\frac{1}{2}k}{0}{- \frac{1}{k} \leq x \leq \frac{1}{k}}
	\ifreport
		\\ 
		\hspace{3em}	
	\fi
	\text{\ \ with \ \ }
	H_k(x) = 	\begin{cases}
					0 								& \text{ if } x < -\frac{1}{k} \\
					\frac{1}{2} + \frac{1}{2}k x	& \text{ if } - \frac{1}{k} \leq x \leq \frac{1}{k} \\
					1								& \text{ if } x > \frac{1}{k}
				\end{cases}
\end{aligneq}
where $H_k(x)$ represents a continuous approximation of the Heaviside function, for large enough $k$:
This formulation satisfies the conditions in \cref{eq:dirac_properties}, and it represents the limit version of the result in \cref{eq:heaviside_dirac_relation}.
Since $h$ is the smallest interval of time we have, the impulse is approximated as $\delta(t-\tau_d) \approx \dert{H_{1/h}}(t-\tau_d)$, where $\dert{H_{{1/h}}}$ is approximated as in \cref{eq:numerical_approximations}.
The smaller $h$ is, the more accurate the approximations are.

This approach works remarkably well for the bouncing ball model, producing an \emph{exactly equal} trace shown in \cref{fig:sample_bouncing_ball_trace} b).
Numerical approximations different than those in \cref{eq:numerical_approximations} may make the traces differ.
One immediate effect is that the values plotted for the impulsive force change by several orders of magnitude. 
Due to the numerical approximation for the derivative, the smaller $h$ is, the larger the impulse will be.
While smaller $h$ increases the accuracy, it may cause overflows due to the limits of machine precision.
Note that the time at which the ball touches the floor still needs to be accurately located \shortcite{Zhang2008} by adjusting $h$ and our simulator does this (both the symbolic and the numerical versions) which explains why \cref{fig:sample_bouncing_ball_trace} b) has different values for the impulses.

Now consider the following artificial signal and its distributional derivatives up to order $n$:
\ifreport
	\small
\fi
\begin{aligneq}\label{eq:derivatives_artificial_signal}
S(t) = H(t-\tau_d) \ ; \ 
\dert{S}(t) = \delta(t-\tau_d) \ ; \ 
\dernt{S}{2}(t) = \dert{\delta}(t-\tau_d) \ ; \ 
\ldots   \ ; \ 
\dernt{S}{n}(t) = \dernt{\delta}{n-1}(t-\tau_d) \\
\end{aligneq}
\ifreport
	\normalsize
\fi
We compare this signal with the approximation $\tilde{S}(t) \approx H_{1/h}(t-\tau_d)$ from \cref{eq:alternative_construct_delta}, differentiated according to the approximation in \cref{eq:numerical_approximations}. 
For example, the derivative $\dert{H}_{1/h}(t-\tau_d)$ is approximated as:
$$
\dert{H}_{1/h}(\tau_d) \approx \frac{H_{1/h}(\tau_d) - H_{1/h}(\tau_d - h)}{h} \approx \frac{1}{h}
$$

\cref{tab:experiment_delta} gives the values of the these approximations around time $\tau_d$.
The table shows that $(n+1)$ time steps are necessary in the numerical simulator to represent the information that is symbolically stored in an impulse derivative of order $n$.
This shifts the numerical solution by $n \times h$ time units, which may cause significant errors.
Intuitively, the symbolic impulse derivatives represent a compact version of what the numerical equivalent would do across the multiple infinitesimal time steps that are abstracted away.

\begin{table}[htb]
\ifreport
	\small
\fi
\caption{Approximated derivatives of $S(t)$ (\cref{eq:derivatives_artificial_signal}).}\label{tab:experiment_delta}
\centering
\begin{tabular}{lrrrrrrr}
\hline
Time		& $S(t)$& $\dert{S}(t)$	& $\dernt{S}{2}(t)$	& $\dernt{S}{3}(t)$	& \ldots & $\dernt{S}{n-1}(t)$	& $\dernt{S}{n}(t)$			\\ \hline
$\tau_d-h$		& $0$	& $0$			& $0$				& $0$				& \ldots & $0$					& $0$						\\
$\tau_d$		& $1$	& $1/h$			& $1/h^2$			& $1/h^3$			& \ldots & $1/h^n$			& $1/h^n$					\\
$\tau_d+h$		& $1$	& $0$			& $-1/h^2$			& $-2/h^3$			& \ldots & $-(n-2)/h^n$ 		& $-(n-1)/h^n$				\\
$\tau_d+2h$	& $1$	& $0$			& $0$				& $1/h^3$			& \ldots & $\binom{n-2}{2}/h^n$ & $\binom{n-1}{2}/h^n$		\\
$\tau_d+3h$	& $1$	& $0$			& $0$				& $0$				& \ldots & $-\binom{n-2}{3}/h^n$	& $-\binom{n-1}{3}/h^n$		\\
\ldots		& \ldots& \ldots		& \ldots			& \ldots			& \ldots & \ldots 				& \ldots					\\
$\tau_d+(n-1)h$& $1$	& $0$			& $0$				& $0$				& \ldots & $0$ 					& $(-1)^{n-1}\binom{n-1}{n-1}/h^n$	\\
$\tau_d+nh$	& $1$	& $0$			& $0$				& $0$				& \ldots & $0$ 					& $0$						\\
\hline
\end{tabular}
\end{table}

In the approximation of $\dernt{S}{n}(t)$, the table, resembling Pascal's triangle, yields the maximum magnitude of the values computed:
\begin{aligneq}\label{eq:max_magnitude}
\binom{k-1}{k-1}/h^k   \text{\ \  where \ \ } k = \floor\pargroup{\frac{n}{2}}
\end{aligneq}
\cref{eq:max_magnitude} can be used to get a rough estimate of the maximum magnitude of values computed in a simulation of an impulsive differential equation.
For example, if the maximum impulse derivative order is $n$, and the maximum amplitude of any discontinuity in the simulation is $D$, then the maximum magnitude is $D \binom{k-1}{k-1}/h^k$.
This result follows from \cref{eq:max_magnitude} and the properties of distributions, described in \cref{sec:distributions} and \cref{sec:hybridcbds:semantics}.

To summarize, for models that contain no impulse derivatives, the numerical approximation is equivalent to the symbolic manipulation. 
\ifnotreport
	This result is shown in more detail in the technical report \shortcite{Gomes2016c}.
\fi
For models that have impulse derivatives, the symbolic approach is more accurate than the numerical one, due to the delay introduced, which is proportional to the highest derivative order of the impulses.

\section{Conclusion}
\label{sec:general_cbds_conclusion}

This work has explored the use of Dirac delta impulses for the modelling and simulation of hybrid systems represented by impulsive differential equations.
\cref{sec:general_cbds_discussion} compares a simulator that approximates impulses numerically with a simulator that encodes the impulses explicitly in the signal.
We conclude that for models that have no impulse derivatives, both approaches are exactly the same, for the approximations in \cref{eq:numerical_approximations}.
The models that \emph{we} have used \emph{do not} have impulse derivatives and we failed at finding models of physical systems that include impulse derivatives.
The numerical approach described in \cref{sec:general_cbds_discussion} is also more efficient, due to a simpler encoding.
The encoding is a stream of real numbers, as opposed to \cref{eq:generalized_signal_representation}.

For models that have impulse derivatives however, the numerical approach introduces a delay in the signal, proportional to the highest derivative order of the impulses, which causes inaccuracies.
Furthermore, if the models are to be run on an embedded platform, where overflows pose a more significant risk, then the symbolic approach may be used. 
\cref{eq:max_magnitude} provides a rough estimate of the magnitudes involved when using the numerical approach.
As can be seen, the smaller the simulation step size $h$ is, the larger \cref{eq:max_magnitude} will be.
This poses an interesting challenge because in general $h$ has to be small in order to get accurate results.
For the symbolic approach, we hypothesize that it may require less bits to perform the same computations on an embedded platform, because the magnitudes involved may not be dominated by the impulse coefficients.
Finally, the symbolic approach has the benefit of being able to identify when operations that are not defined in theory, are being computed (e.g., an inversion of an impulse). 
As \cite{Nilsson2003} points out, this is not fool proof because the signals are computed point-wise but, compared to the numerical approach, it is an improvement.

\ifnotreport
	\vspace{-10pt}
\fi
\section*{Acknowledgment}
We are thankful to the anonymous reviewers for the corrections and comments provided.
This research was partially supported by Flanders Make vzw, the strategic research centre for the manufacturing industry,
and PhD fellowship grants from (1) the Agency for Innovation by Science and Technology in Flanders (IWT), and (2) Research Foundation - Flanders (FWO).

\ifreport
	\bibliographystyle{plain}
	\bibliography{./library}
\else
	\vspace{-10pt}
	\bibliographystyle{scsproc}
	{\small
	\bibliography{./library}

\begin{thebibliography}{10}

\bibitem{Cellier1991}
Fran{\c{c}}ois~Edouard Cellier.
\newblock {\em {Continuous system modeling}}.
\newblock Springer Science {\&} Business Media, 1991.

\bibitem{Cellier2006}
Fran{\c{c}}ois~Edouard Cellier and Ernesto Kofman.
\newblock {\em {Continuous System Simulation}}.
\newblock Springer Science {\&} Business Media, 2006.

\bibitem{Dirac1981}
Paul Adrien~Maurice Dirac.
\newblock {\em {The principles of quantum mechanics}}.
\newblock Number~27. Oxford university press, 1981.

\bibitem{Friedlander1998}
Friedrich~Gerard Friedlander and Mark~Suresh Joshi.
\newblock {\em {Introduction to the Theory of Distributions}}.
\newblock Cambridge University Press, 1998.

\bibitem{Gomes2016a}
Cl{\'{a}}udio Gomes, Joachim Denil, and Hans Vangheluwe.
\newblock {Causal-Block Diagrams}.
\newblock Technical report, 2016.

\bibitem{Horvath1970}
John Horvath.
\newblock {An Introduction to Distributions}.
\newblock {\em The American Mathematical Monthly}, 77(3):227, mar 1970.

\bibitem{Lakshmikantham1989}
Vangipuram Lakshmikantham, Drumi~D Bainov, and Pavel~S Simeonov.
\newblock {\em {Theory of impulsive differential equations}}, volume~6.
\newblock World scientific, 1989.

\bibitem{Lee2015}
Edward~A Lee, Mehrdad Niknami, Thierry~S Nouidui, and Michael Wetter.
\newblock {Modeling and Simulating Cyber-physical Systems Using CyPhySim}.
\newblock In {\em Proceedings of the 12th International Conference on Embedded
  Software}, EMSOFT '15, pages 115--124, Piscataway, NJ, USA, 2015. IEEE Press.

\bibitem{Mosterman1998a}
Pieter~J. Mosterman and Gautam Biswas.
\newblock {A theory of discontinuities in physical system models}.
\newblock {\em Journal of the Franklin Institute}, 335(3):401--439, apr 1998.

\bibitem{Mustafiz2016a}
Sadaf Mustafiz, Cl{\'{a}}udio Gomes, Bruno Barroca, and Hans Vangheluwe.
\newblock {Modular Design of Hybrid Languages by Explicit Modeling of Semantic
  Adaptation}.
\newblock In {\em Proceedings of the Symposium on Theory of Modeling {\&}
  Simulation: DEVS Integrative M{\&}S Symposium}, DEVS '16, pages 29:1----29:8,
  San Diego, CA, USA, 2016.

\bibitem{Nilsson2003}
Henrik Nilsson.
\newblock {Functional automatic differentiation with dirac impulses}.
\newblock {\em ACM SIGPLAN Notices}, 38(9):153--164, sep 2003.

\bibitem{Posse2002}
Ernesto Posse, Juan de~Lara, and Hans Vangheluwe.
\newblock {Processing causal block diagrams with graphgrammars in atom3}.
\newblock In {\em European Joint Conference on Theory and Practice of Software
  (ETAPS), Workshop on Applied Graph Transformation (AGT)}, pages 23--34, 2002.

\bibitem{Stewart2000}
D~Stewart.
\newblock {Rigid-Body Dynamics with Friction and Impact}.
\newblock {\em SIAM Review}, 42(1):3--39, jan 2000.

\bibitem{Strichartz2003}
Robert~S Strichartz.
\newblock {\em {A guide to distribution theory and Fourier transforms}}.
\newblock World Scientific, 2003.

\bibitem{Vangheluwe2014}
Hans Vangheluwe, Joachim Denil, Sadaf Mustafiz, Daniel Riegelhaupt, and Simon
  {Van Mierlo}.
\newblock {Explicit Modelling of a CBD Experimentation Environment}.
\newblock In {\em Proceedings of the Symposium on Theory of Modeling {\&}
  Simulation - DEVS Integrative}, DEVS '14, San Diego, CA, USA, 2014. Society
  for Computer Simulation International.

\bibitem{Zhang2008}
Fu~Zhang, Murali Yeddanapudi, and Pieter Mosterman.
\newblock {Zero-crossing location and detection algorithms for hybrid system
  simulation}.
\newblock In {\em IFAC World Congress}, pages 7967--7972, 2008.

\end{thebibliography}
	}
\fi

\end{document}